\documentclass[11pt,a4paper,reqno]{amsart}
\usepackage[english]{babel}
\usepackage[T1]{fontenc}
\usepackage{palatino}
\usepackage{amsmath}
\usepackage{amssymb}
\usepackage{amsthm}
\usepackage{amsfonts}
\usepackage[utf8]{inputenc}
\usepackage{graphicx}
\usepackage{color}
\usepackage{mathtools}

\usepackage[colorlinks = true, citecolor = blue]{hyperref}
\title{Infinite geodesics of sub-Finsler distances in Heisenberg groups}
\author[Z.\ M.\ Balogh]{Zolt\'{a}n M.\ Balogh}
\address{Department of Mathematics and Statistics\\ University of Bern \\ Sidlerstrasse 5, 3012 Bern, Switzerland}
\email{zoltan.balogh@math.unibe.ch}
\author[A. Calogero]{Andrea Calogero}
\address{Dipartimento di Matematica e Applicazioni\\ Universit\'a di Milano Bicocca \\ Via Cozzi 53, 20125 Milano, Italy}
\email{andrea.calogero@unimib.it}

\date{\today}
\subjclass[2010]{53C17, 22E25, 49K15}
\thanks{Z.M.B.\ was supported by the Swiss National Science Foundation through the project 165507 `Geometric Analysis of Sub-Riemannian Spaces'.}
\keywords{Heisenberg group, homogeneous norms, geodesics, optimal control}
\newcommand{\R}{\mathbb{R}}

\newcommand{\E}{\mathbb{E}}

\def\Barint_#1{\mathchoice
          {\mathop{\vrule width 6pt height 3 pt depth -2.5pt
                  \kern -8pt \intop}\nolimits_{#1}}%
          {\mathop{\vrule width 5pt height 3 pt depth -2.6pt
                  \kern -6pt \intop}\nolimits_{#1}}%
          {\mathop{\vrule width 5pt height 3 pt depth -2.6pt
                  \kern -6pt \intop}\nolimits_{#1}}%
          {\mathop{\vrule width 5pt height 3 pt depth -2.6pt
                  \kern -6pt \intop}\nolimits_{#1}}}
\def\QED{\hfill $\Box$\par\smallskip\noindent}

\def\scatola{\lower5pt\hbox{\vbox{\hrule\hbox{\vrule\kern2pt\vbox%
{\kern5pt\hbox{\mathsurround=0pt }\kern2pt}\kern4pt\vrule}\hrule}}\
} 

\def\de{{\rm \, d}}
\def\H{{\mathbb H}}
\def\Hn{{\mathbb H^n}}

\def\R{{\mathbb{R}}}
\def\Rn{{\mathbb{R}^n}}

\numberwithin{equation}{section}

\theoremstyle{plain}
\newtheorem{thm}[equation]{Theorem}

\newtheorem{proposition}[equation]{Proposition}

\newtheorem{example}{Example}[section]
\theoremstyle{definition}

\newtheorem{definition}[equation]{Definition}

\theoremstyle{remark}
\newtheorem{remark}[equation]{Remark}

\begin{document}

\begin{abstract}
 We consider Heisenberg groups equipped with a sub-Finsler metric. Using methods of optimal control theory we prove that in this geometric setting the infinite geodesics are horizontal lines under the assumption that  the sub-Finsler metric is defined by a strictly convex norm. This answers a question posed in \cite{BalFasSob2017} and has applications in the characterisation of isometric embeddings into Heisenberg groups.   \end{abstract}
\maketitle

\section{Introduction}

In the recent  paper  \cite{BalFasSob2017} the problem of classification of  isometric embeddings of Heisenberg groups $\mathbb{H}^m$ into $\mathbb{H}^n$ for $m\leq n$ has been considered.
Here both groups  $\mathbb{H}^m$ and $\mathbb{H}^n$ were endowed with a  homogeneous distance.  By such a distance, we mean a left-invariant metric induced by a gauge function which is homogeneous with respect to
a one-parameter family of `Heisenberg dilations' adapted to the stratification of the underlying Lie algebra. These type of geometric structures are interesting objects  and  have been in the focus of several recent studies as shown in  \cite{ABB,BBDS2017,  BelRis1996, DLR2017, LeRi2017, HaZi2015, Mon2000}. In the above referenced papers the more local aspects of the geometry were mainly considered. Our purpose is to concentrate on more global geometrical aspects as we study the behaviour of infinite geodesics in sub-Finsler distances.

The main motivation of our work is the following result of \cite{BalFasSob2017}: consider
 $\mathbb{H}^m$ and $\mathbb{H}^n$, $n\geq m$, to be equipped with left-invariant homogeneous distances $d$ and $d'$, respectively. Assuming that every infinite geodesic in
 the target space $(\mathbb{H}^n,d')$ is a line, then every isometric embedding $f:(\mathbb{H}^m,d) \to (\mathbb{H}^n,d')$ is the composition of a left translation and a homogeneous homomorphism.

 \medskip

 The above result raises the natural question of characterising homogeneous distances $d$ defined on $\mathbb{H}^n$ such that the space $(\mathbb{H}^n,d)$ has the property that every infinite geodesic in $(\mathbb{H}^n,d)$ is a line. In \cite{BalFasSob2017} this property was called the geodesic linearity property or $GLP$ of the space $(\mathbb{H}^n,d)$.
 It was conjectured in \cite{BalFasSob2017} that $(\mathbb{H}^n,d)$ satisfies $GLP$ if and only if the
 underlying Euclidean norm $N$ associated to $d$ is strictly convex. In \cite{BalFasSob2017} this conjecture has been verified for several concrete examples.  The purpose of this paper is to prove this conjecture for general homogeneous distances on $\mathbb{H}^n$. In order to formulate our result we need to fix notation and recall some preliminaries.


The \emph{$n$-th Heisenberg group}  $\mathbb{H}^{n}$ is the set $\mathbb{R}^{2n}\times\mathbb{R}$ equipped with the
multiplication
\begin{displaymath}
g\ast g'= (z,t)\ast(z',t'):=(z+z',t+t'+2\langle z,J_n z'\rangle), \text{ where } J_n = \begin{pmatrix}0&-E_n\\E_n&0\end{pmatrix}\in\mathbb{R}^{2n\times2n},
\end{displaymath}
is the standard symplectic matrix on $\mathbb{R}^{2n}$ and $E_n$ denotes the $(n\times n)$ unit matrix. Sometimes it is convenient to write in coordinates
\begin{displaymath}
z=(x_1,\ldots,x_n,y_1,\ldots,y_n).
\end{displaymath}
It can be easily verified
that $(\mathbb{H}^{n},\ast)$ satisfies all properties of a group
with neutral element $e:=(0,0)$ and inverse $(z,t)^{-1}:=(-z,-t)$.  Denoting the nonlinear term $\langle z,J_n z'\rangle$ by $\omega_{n} (z,z')$, we remark that this expression defines a skew-symmetric bilinear form on $\mathbb{R}^{2n}$, and that two elements $(z,t)$ and $(z',t')$ in $\mathbb{H}^{n}$ commute if and only if the term $\omega_n(z,z')$ is zero. Since this does not hold for all elements in $\mathbb{H}^{n}$, it turns out that the Heisenberg group is  non-abelian.
For $\lambda>0$, the map $\delta_{\lambda}:\mathbb{H}^{n}\to\mathbb{H}^{n},\ (z,t)\mapsto(\lambda z,\lambda^{2}t)$ is called \emph{$\lambda$-dilation}.
It can be easily verified that any $\lambda$-dilation defines a group isomorphism with inverse $\delta_{\lambda^{-1}}$. It plays an analogous role as the usual scalar multiplication in $\mathbb{R}^n$.

 We say that a  \emph{norm} on $\Hn$ is a map $\hat{N}:\Hn \to\mathbb{R}_{\geq0}$ that satisfies
\begin{align}
&i)\ \hat{N}(g)=0\Leftrightarrow g=e,\ \forall g\in \Hn,\nonumber\\
&ii)\ \hat{N}(g^{-1})=\hat N(g),\ \forall g\in \Hn,\nonumber\\
&iii)\ \hat{N}(g*g')\leq \hat{N}(g)+\hat{N}(g'),\ \forall g,g'\in \Hn.\nonumber
\end{align}

 A metric $d:\Hn\times \Hn \to\mathbb{R}_{\geq0}$
is called \emph{left-invariant}, if for every $g_{o}\in \Hn$, the map $L_{g_{o}}:(\Hn,d)\to(\Hn,d),\ g\mapsto g_{o}\ast g$ is an isometry, that is, $d(g_{o}\ast g,g_{o}\ast g')=d(g,g')$, for all $g,g'\in \Hn$.

Every norm $\hat{N}:\Hn\to\mathbb{R}_{\geq0}$ induces
a left-invariant metric $d_{\hat{N}}:\Hn\times \Hn\to\mathbb{R}_{\geq0}$,
and vice versa. More precisely, we can establish the following bijection
\begin{align}
\{\hat{N}:\Hn\to\mathbb{R}_{\geq0}:\;\hat N\ \text{is a norm}\}&\to\{d:\Hn \times \Hn\to\mathbb{R}_{\geq0}:\;d\ \text{is a left-invariant metric}\}\nonumber\\
N&\mapsto d_{\hat N}:\Hn\times \Hn\to\mathbb{R}_{\geq0},\ (g,g')\mapsto \hat N(g^{-1}*g'),\nonumber
\end{align}
\begin{align}
\{d:\Hn\times \Hn\to\mathbb{R}_{\geq0}:\;d\ \text{is a left-invariant metric}\}&\to\{\hat N:\Hn\to\mathbb{R}_{\geq0}:\;N\ \text{is a norm}\}\nonumber\\
d&\mapsto \hat N_{d}:\Hn\to\mathbb{R}_{\geq0},\ g\mapsto d(g,e).\nonumber
\end{align}

A norm $\hat N:\mathbb{H}^{n}\to\mathbb{R}_{\geq0}$ on the Heisenberg group is called \emph{homogeneous} if
\begin{align}
\ \hat N(\delta_{\lambda}(g))=\lambda \hat N(g),\ \text{for all }\lambda>0,\ \text{for all } g\in\mathbb{H}^{n}.\nonumber
\end{align}

It is easy to see that a norm $\hat N$ on $\mathbb{H}^{n}$ is homogeneous if and only its associated left-invariant metric is homogeneous in the sense that $d_{\hat N}(\delta_{\lambda}(g),\delta_{\lambda}(g'))=\lambda d_{\hat N }(g,g')$. 
Every left-invariant distance on $\mathbb{H}^n$ induced by a homogeneous norm is a homogeneous distance.
From now on, we will use the expression "\emph{homogeneous distance on $\mathbb{H}^{n}$}" to talk about the left-invariant metric induced by a homogeneous norm. It follows from that the topology induced by any homogeneous distance on $\mathbb{H}^{n}$ coincides with the Euclidean topology on $\mathbb{R}^{2n+1}$, and  that any homogeneous norm is  continuous with respect to the Euclidean topologies of $\mathbb{R}^{2n+1}$ and $\mathbb{R}$.
In particular, we note that any two homogeneous distances on $\mathbb{H}^{n}$ induce the same topology.

\medskip

 On the other hand, the metric structure induced by a homogeneous norm $\hat N$ on $\mathbb{H}^n$ is very different from $\mathbb{R}^{2n+1}$ endowed with the Euclidean distance $d_{eucl}$. The two distances $d_{\hat N}$ and $d_{eucl}$ are not bi-Lipschitz equivalent for any choice of homogeneous norm $\hat N$ on $\mathbb{H}^n$. However, we can associate to $\hat N$ a norm $N$ by restricricting it to $\R^{2n}$ by:
 \begin{equation} \label{norm-euclidean}
 N(z):= \hat N(z,0) , \qquad \text{ for} \ z\in \R^{2n}.
 \end{equation}
 Indeed, one can easily check (see Proposition 2.8 \cite{BalFasSob2017}) that $N: \R^{2n} \to \R$ satisfies the axioms of a norm on the  vector space $\R^{2n}$.
 The main result of the paper is the following:

 \begin{thm} \label{main}Let $(\mathbb{H}^n, d)$ be the Heisenberg group endowed with the homogeneous distance $d$. Denote by $\hat N: \Hn \to \R_{\ge 0}$ the associated homogeneous norm such that $d= d_{\hat N}$ and consider the associated norm $N: \R^{2n}\to \R_{\ge 0}$ defined by \eqref{norm-euclidean}. Then the space  $(\mathbb{H}^n, d)$ has the GLP if and only if $N$ is strictly convex.
  \end{thm}

 Let us recall that by definition, the norm $N:\R^{2n}\to \R_{\ge 0}$ is strictly convex iff its closed unit ball $B_N$ (see Definition \ref{BN} below) is a strictly convex set, that is, if $z_1,z_2 \in B_N$ and $\alpha \in (0,1)$ then $N(\alpha z_1 + (1-\alpha) z_2) < 1$.

Recall also that by  a \emph{geodesic} $\gamma:I \to (\Hn,d)$, we mean an isometric embedding of $I=[a,b]$ or $I=\mathbb{R}$ into $(\Hn,d)$, that is,
\begin{equation}\label{def geodesic}
d(\gamma(s),\gamma(s'))=|s-s'|,\quad\text{for all }s,s'\in I.
\end{equation}
If we have  $I=\mathbb{R}$ in the above definition, we say that $\gamma$ is an \emph{infinite geodesic}.

\medskip

To comment on the statement of Theorem \ref{main} let us mention that in the setting of normed spaces the strict convexity of the norm is equivalent to the fact that {\bf all} (and not just the infinite geodesics) are straight line segments or lines (see eg. \cite{Pap2005}). This fact is definitely false in the sub-Finsler setting of the Heisenberg group. Let us recall (\cite{Mon2000}) that in the standard sub-Riemannian metric the underlying norm is the Euclidean one, that is strictly convex. On the other hand there are a multitude of finite geodesics that are not line segments.

\medskip

In order to prove Theorem \ref{main} we shall consider the  the \emph{sub-Finsler distance} associated to $N$. The definition of the sub-Finsler distance is based on the notion of horizontal curves. A \emph{horizontal curve} in $\mathbb{H}^n$ is an absolutely continuous curve $\gamma:[a,b] \to \mathbb{R}^{2n+1}$ with the property that
\begin{displaymath}
\dot{\gamma}(s) \in H_{\gamma(s)}, \quad\text{for almost every }s\in [a,b],
\end{displaymath}
where for $g\in \mathbb{H}^n$, we set
\begin{displaymath}
H_g:= \mathrm{span}\left\{X_{1,g},\ldots,X_{n,g},Y_{1,g}\ldots,Y_{n,g} \right\}.
\end{displaymath}
Here $X_i$ and $Y_i$, $i=1,\ldots,n$, are the left-invariant vector fields (with respect to $\ast$) which at the origin agree with the standard basis vectors: $X_{i,e}=e_i$ and $Y_{i,e}= e_{n+i}$.
These left-invariant vector fields can be written as first order
differential operators as follows
\begin{equation} \label{vector fields} \left.
  \begin{array}{lll}
     X_j=\partial_{x_j}+2y_j \partial_t,\qquad j=1,...,n,\\
     Y_j=\partial_{y_j}-2x_j\partial_t,\qquad j=1,...,n.
  \end{array}
\right.
\end{equation}

Denoting the $(2n+1)$ components of an absolutely continuous curve $\gamma:[a,b]\to \mathbb{H}^n$ by $\gamma_i$, $i=1,\ldots,2n+1$, it follows that $\gamma$ is horizontal if and only if
\begin{equation}\label{eq:horiz_comp}
\dot{\gamma}_{2n+1}(s) = 2 \sum_{i=1}^n \left(\dot{\gamma}_i(s) \gamma_{n+i}(s)-\dot{\gamma}_{n+i}(s)\gamma_i(s)\right),\quad\text{for almost every }s\in [a,b].
\end{equation}
It is well known that a horizontal curve $\gamma:[a,b]\to \mathbb{H}^n$ is rectifiable and admits a Lipschitz parametrization (see for instance \cite[Proposition 1.1]{HaZi2015} for a proof and  note that this statement holds for any homogeneous norm on $\mathbb{H}^n$).
 In converse direction, every rectifiable curve admits a $1$-Lipschitz parametrization and this parametrization is horizontal, see \cite{Pa1989}.

Given a norm $N:\R^{2n}\to \R_{\ge 0},$ the \emph{sub-Finsler distance} associated to $N$ on $\mathbb{H}^n$ is given by
\begin{equation}\label{sub-finsler}
d_{SF}(g,g'):= \inf_{\gamma} \int_a^b N(\dot{\gamma}_I(s)) \;\mathrm{d}s,
\end{equation}
where the infimum is taken over all horizontal curves $\gamma=(\gamma_{I},\gamma_{2n+1}): [a,b] \to \mathbb{H}^n$ with $\gamma(a)=g$ and $\gamma(b)=g'$, $\gamma_I= (\gamma_1, \ldots, \gamma_{2n})$. Here, and in the remaining part of the paper we shall refer to $\gamma_I$ as the projection of the horizontal curve $\gamma$ and if $\gamma_I$ is given the curve $\gamma= (\gamma_I, \gamma_{2n+1})$ satisfying \eqref{eq:horiz_comp} is referred to as the horizontal lift of $\gamma_I$. It is important to mention that given $g,g' \in \H^n$ and
$a< b \in \R $ with $b-a= d_{SF}(g,g')$ then the above minimisation problem \eqref{sub-finsler} has always a solution  and horizontal curves solving this problem are (up to a  reparametrization) geodesics in the space $(\Hn, d_{SR})$.

\medskip
Coming back to Theorem \ref{main} let us mention that by Proposition 3.14, \cite{BalFasSob2017} we alread have one direction of the statement, namely that $GLP$ of $(\H^n, d_{\hat N})$ implies that $N$ is strictly convex. Our purpose here is to prove the other implication that is substantially more difficult. Namely, we need to show that if $N$ is strictly convex, then every infinite geodesic in $(\H_n, d_{\hat N})$ is a straight line. Notice first that according to Proposition 2.19 in \cite{BalFasSob2017} if $\gamma: \R \to \mathbb{H}^n$ is an infinite  geodesic with respect to $d_{\hat N}$, then it is also an infinite geodesic with respect to $d_{SF}$.
Therefore we shall only need to consider infinite geodesics associated to $d_{SF}$.
In order to study properties of geodesics related to the sub-Finsler distance \eqref{sub-finsler} we reformulate this minimisation problem as a problem of optimal control theory.

For the sake of simplicity we assume that $a= 0, b= T$, the starting point $\gamma(0)=e$ is the neutral element of $\H^n$ and the final point is
$\gamma(T)=g\not = e$. Then the length-minimising property of the geodesic from  \eqref{sub-finsler} will be equivalent to the following optimal control problem with fixed time and fixed end-point:
\begin{equation}\label{control problem}
\left\{
  \begin{array}{ll}
     \displaystyle \inf_v \int_0^T N(v(s))\de s\\
\dot x_i= -v_i,\qquad i=1,\ldots, n\\
\dot y_i= -v_{n+i},\qquad i=1,\ldots, n\\
\displaystyle \dot{t}=-2\sum_{i=1}^n\left(y_{i}v_i-x_iv_{n+i}\right)\\
\gamma(0)=e\\
\gamma(T)=g
  \end{array}
\right. .
 \end{equation}
 Notice, that in the above formulation, the curve   $\gamma=(x_1,\ldots,x_n,y_1,\ldots,y_n,t):[0,T]\to \Hn$ is automatically a horizontal curve, and the control function $-v:[0,T]\to \Rn$ is in fact the horizontal velocity of $\gamma$, $\dot \gamma_I =- v$.
 It is clear that the two minimisation problems \eqref{sub-finsler} and \eqref{control problem} are equivalent to each other.

Our approach is to apply Pontryagin's Maximum Principle in order to obtain useful information about solutions of \eqref{control problem}. We should mention here that the study of sub-Finsler metrics by this approach has also been undertaken by \cite{BBDS2017}. However, in the aforementioned  work, properties of infinite geodesics were not addressed. In fact, the main technical difficulty of our paper is precisely to pass from finite to infinite geodesics because of the possible ambiguity of the multiplier arising in Pontryagin's theorem.

The paper is organized as follows: in Section 2 we review some necessary background on convex analysis that is going to be used in the sequel. In Section 3 we apply Pontryagin's theorem to our situation. Section 4 contains the proof of our main result and Section 5 is for final remarks and examples.


\section{Preliminary results from convex analysis}\label{section concex analysis}

Our goal is to study a class of norms on $\Hn$ that are generated by convex norms
coming from $\R^{2n}$. In this section we collect some basic results from the Euclidean convex analysis to be used in the sequel.

Let us start with a general norm $N:\R^{2n}\to \R_{\ge 0}$ and let us denote by $B_N$ its unit ball, i.e. the set
\begin{equation}\label{BN}
B_N:=\{z\in\R^{2n}:\ N(z)\le 1\}.
\end{equation}
It is clear that $B_N$ is compact, convex with $0$ in its interior.
Since for every $\lambda\in [0,1]$ and $z,\ z'$ in $\R^{2n}$
$$
N((1-\lambda)z+\lambda z')\le N((1-\lambda)z)+N(\lambda z')=(1-\lambda)N(z)+\lambda N(z'),
$$
every norm is a convex function in $\R^{2n}$. We define the \emph{dual norm} $N_*$ of $N$ as usual by
$$ N_*(p) = \sup \{ p\cdot z:\  N(z) = 1,\ z\in \R^{2n} \}.$$
The  unitary ball $B_{N_*}$ of $N_*$ is the polar set $B_N^\circ$ of $B_N$, i.e.
$$
B_N^\circ:=\{p\in\R^{2n}:\ p\cdot z\le 1,\ z\in B_N\}.$$

 Given a function $f:\R^{2n}\to \R$, the  \emph{subdifferential} is a set valued mapping  $\partial f: \R^{2n} \rightrightarrows \R^{2n}$ defined by
$$ \partial f(z)= \{ p \in \R^{2n}: f(z') \geq f(z) + p\cdot(z'-z), \ \forall\, z' \in \R^{2n} \}.$$
A well-known result due to Rockafellar characterizes the convexity of a function via the subdifferential; more precisely,
$f:\R^{2n}\to\R$ is convex if and only if
$\partial f(z) \neq \emptyset $ for all $z \in \R^{2n}$ (see for instance \cite{RoWe2004}).     Moreover, for  a convex function $f$,
the subdifferential set valued mapping $\partial f$ is upper semicontinuous, and it is compact and convex valued
(see \cite[Proposition 2.1.5 and Proposition 2.2.7]{Cl1983} and also \cite{AlAmCa1992}).
These properties of the map $\partial f$ will play a fundamental role in the case of non-smooth norm (see Proposition \ref{proj-estimate} below).

The \emph{Legendre transform} of $f$ is $f^*:\R^{2n}\to (-\infty,+\infty]$ defined by
$$
f^*(p)=\sup_{z\in\R^{2n}}\left(p\cdot z -f(z)\right)
$$
Then $f^*$ is always a convex function;  if $f$ is superlinear, then $f^*$ is real valued.
For every convex function $f$, the Legendre transform $f^*$ is related to the subdifferential
via the following equality (see for example \cite[Proposition 11.3]{RoWe2004})
\begin{equation}\label{Fenchel}
 p\in\partial f(z)\  \Longleftrightarrow \
z\in\partial f^*(p)\ \Longleftrightarrow \
 f(z)+f^*(p)=z\cdot p.
\end{equation}
In all that follows, given a norm $N$ we associate the function $F_N:\R^{2n}\to\R_{\ge 0}$ by
\begin{equation}\label{FN}
F_N(z)=\frac{1}{2}(N(z))^2.
\end{equation}
As we shall see, $F_N$ will play a crucial role in this paper.

The previous notions and properties are closely intertwined with the homogeneity and the convexity of the norm $N$ as follows:
\begin{proposition}\label{Propo first}
Let $N$ be a norm in $\R^{2n}$. The following properties hold:
\begin{itemize}
\item[i.] for all $z\in\R^{2n}$
\begin{equation}\label{subgradient}
p\in\partial N(z)\qquad\Longleftrightarrow \qquad p\in\partial N(\alpha z), \quad \forall \alpha>0;
\end{equation}

\item[ii.]  for every $z\in\R^{2n},$ with $z\not=0$, we have
\begin{equation}\label{Eulero}
  p\in\partial N(z)\qquad  \Longleftrightarrow \qquad
 N(z)=z\cdot p.
\end{equation}
\item[iii.] $F_N$ is a convex function and hence $\partial F_N(z)\not=\emptyset,$ for every $z\in\R^{2n}$;
moreover, for every $z\not=0$
\begin{equation}\label{chain rule}
\partial F_N(z)=N(z)\partial N(z);
\end{equation}

\item[iv.] The Legendre transform $(F_N)^*$ of $F_N$ is denoted by  $F_N^*$ and satisfies the relation $F^*_N=F_{N_*},$ i.e.
$$
F_N^*(p)=\frac{1}{2}(N_*(p))^2.
$$
\end{itemize}

\end{proposition}
\noindent\textbf{Proof:} The proof of \eqref{subgradient} is an easy exercise.

Let us pass to ii: based on \eqref{subgradient} we can assume that $N(z)=1$. Thus we have to prove
that $p\in \partial N(z)$ if and only if $p\cdot z =1$.

According to \eqref{Fenchel} we have that $p\in \partial N(z)$ if and only if $N(z)+ N^*(p) = p\cdot z$. In our case this relation reads as $1+ N^*(p) = p\cdot z$.
Let us consider the  Legendre trasform $N^*$ of $N$.
Since
\begin{eqnarray*}
N^*(p)&=&\sup_{z\in\R^{2n}}\left(p\cdot z -N(z)\right)\\
&=&\sup_{\alpha\ge 0}\alpha\biggl(\sup_{\{z:\, N(z)=1\}}\left(p\cdot z -1\right)\biggr).
\end{eqnarray*}
By the above relation it follows that $N^*$ takes only  values in the set $\{0,\, +\infty\}$. Moreover,  $N^*(p)=0$ if and only if  $p\in B_N^\circ$. This implies that
$$
 p\in\partial N(z)\  \Longleftrightarrow \
z\cdot p= 1,
$$
which proves ii.

Let us prove iii.
For  $z_0$ and $z_1$ in $\R^{2n}$ and $\alpha\in[0,1]$ we have
\begin{eqnarray}
F_N((1-\alpha) z_0+\alpha z_1)
&\le&\frac{1}{2} \left(N((1-\alpha) z_0)+N(\alpha z_1)\right)^2\nonumber\\
&=&\frac{1}{2} \left((1-\alpha)N( z_0)+\alpha N( z_1)\right)^2\label{F1}\\
&\le&\frac{1}{2} \left((1-\alpha)(N( z_0))^2+\alpha (N( z_1))^2\right)\label{F2}
\end{eqnarray}
The inequality \eqref{F1} follows from the homogeneity of the norm $N$ and the fact that $s\mapsto s^2$ is an increasing function in $[0,\infty)$. Inequality \eqref{F2} follows from the convexity of the function $s\mapsto s^2$.

\noindent
Now consider a convex function $f:\R^{2n}\to \R$ and a non decreasing,  convex and regular function $\varphi:I\to \R$, where $I\subset\R$ is an interval such that $f(\R^{2n})\subset I$, then the following version of the chain rule holds (see \cite[Proposition 4.2.5]{BeNeOz2003}):
$$\partial (\varphi\circ f)(z) =  \varphi '(f(z))\partial f(z).$$
Applying this result to our case with $\varphi(s)=\frac{s^2}{2}$ and $f=N$, we have
\eqref{chain rule} and the proof of iii.

For the proof of {iv.} see Proposition 11.21 in \cite{RoWe2004}, taking into account that the function $\theta(s)=\frac{1}{2}s^2$ is convex and $\theta^*=\theta$.
\QED

Let us introduce the following fundamental notion:
\begin{definition}
We say that a norm $N:\R^{2n}\to\R_{\ge 0}$ is \emph{strictly convex} if $z$ and $z'$ in $\R^{2n}\setminus\{0\}$ are such that $N(z+z')=N(z)+N(z')$, then $z'=\alpha z$ for some $\alpha>0$.

\noindent
We say that a norm is \emph{smooth} if at every point $z\in\partial B_N$ the ball $B_N$ has a unique supporting hyperplane.
\end{definition}
As an equivalent notion to require that a norm $N$ is strictly convex iff $\partial B_N$ does not contain line-segments.
Furthermore, the smoothness of $N$ is equivalent to
the fact that, for every $z\not=0$, $N$ is differentiable, or  that $\partial N(z)$ is singleton, i.e.
$$
\partial N(z)=\{\nabla N(z)\},\quad\forall z\not=0.
$$
We will denote by $\|\cdot\|$ the classical Euclidean norm. Finally, let us recall that
a norm $N$ is strictly convex if and only if $N_*$ is smooth.
See for instance Chapter 5 in \cite{Me1998} for more details.

Every norm $N$ is a convex function however,  the strict convexity of $N$ as a norm does not imply its strict convexity as a function: in fact
 $$
 N(\lambda z)=\lambda N(z)+(1-\lambda)N(0),\qquad \forall z\in\R^{2n},\ \lambda\in[0,1].
 $$
 Taking the square of $N$ and by defining the function $F_N$ as in \eqref{FN} will get us around this difficulty and, starting with a strictly convex norm, we obtain a useful strictly convex function. Hence we have the following result:

\begin{proposition}\label{Propo zero bis}
Let $N$ be a strictly convex norm. Then $F_N$ is a strictly convex function and $\partial F_N:\R^{2n}\rightrightarrows \R^{2n}$ is injective, i.e.
$$z\not=z'\quad\Rightarrow\quad \partial F_N(z)\cap\partial F_N(z')=\emptyset.$$
\end{proposition}
\noindent\textbf{Proof:}
Observe first that since  $s\mapsto s^2$ is a strictly increasing function in $[0,\infty)$, in \eqref{F1} we have an equality if and only if
$
N((1-\alpha) z_0+\alpha z_1)=
N((1-\alpha) z_0)+N(\alpha z_1).$ Since $N$ is a strictly convex norm, this implies that there exists $\theta>0$ such that
\begin{equation}\label{F3}
(1-\alpha) z_0=\theta \alpha z_1.
\end{equation}
Since  $s\mapsto s^2$ is a strictly convex function,
in  \eqref{F2} we have an equality if and only if
\begin{equation}\label{F4}
N(z_0)=N(z_1).
\end{equation}
Relations \eqref{F3} and \eqref{F4} give $z_0=z_1$. The injectivity of $\partial F_N$ follows from its strict convexity.  \QED

The strictly convexity of the  function $F_N$ plays a fundamental role in that follows. Indeed, let us recall that
if $f:\R^{2n}\to [0,\infty)$ is a strictly convex function, then it may happen that $f$ does not have a minimum point, i.e. it is not possible to guarantee in general that there exists $z^*$ such that
$$
z^*\in\arg\min_{z\in\R^{2n}}f(z).
$$
However, if such $z^*$ exists, by the strict convexity of $f$, it is unique.

\section{Optimal control approach to geodesics in $\Hn$.}

Let $N:\R^{2n}\to\R_{\ge 0}$ be a norm. The Heisenberg geodesics $\gamma=(x_1,\ldots,x_n,y_1,\ldots,y_n,t):[0,T]\to \Hn$ with respect
 to the norm $N$, joining the points $e, g \neq e \in \Hn$ will be determined by the optimal control problem \eqref{control problem} with fixed final time $T>0$ and fixed end-point $g$
 (which is equivalent to the problem \eqref{sub-finsler}).
 Let us only recall that in the
 formulation  \eqref{control problem}, the curve   $\gamma$ is automatically a horizontal curve, and the control function $-v=[0,T]\to \Rn$ is in fact the horizontal velocity of $\gamma$, $\dot \gamma_I =- v$.

\subsection{Equivalence of optimal control problems related to $N$ and $F_N$.}

If we want to apply the methods of optimal control theory to the problem \eqref{control problem} we shall run into problems due to the fact that the norm $N$ has a linear growth at infinity. Therefore, our purpose is to change the above setting to an equivalent problem related to $F_N$ where the integrand has a quadratic growth at infinity.
The first step in this direction is the following  normalization lemma where the homogeneity of $N$ is crucial.

 \begin{proposition}\label{reparametrization} Let $N:\R^{2n}\to \R_{\ge 0}$ be a norm.
 Let us consider the problem \eqref{control problem}.
Without loss of generality we can reparametrize by arc
length the horizontal curve $\gamma$ such that
$$N(v(s))=K,\qquad s\in[0,T],$$
where $v=-\dot\gamma_I$ and $K$ is a positive constant depending only on $T$ and $\gamma$.
  \end{proposition}

 \noindent \textbf{Proof}:
Consider  a horizontal curve $\gamma : [0, T] \to \Hn$ from $0$ to $g$, with
  $$\ell(\gamma)=\int_0^T N(v(s))\de s$$
  and $v=-\dot\gamma_I$. If the two endpoints are different: $g\not =e$ we can assume without loss of generality that $v(s) \not = 0$ for all $s\in [0,T]$.
  Then we can define an absolute continues homeomorphism $\tau: [0, T]\to[0,T]$ by
  $$
  \tau( s)=\frac{T}{\ell(\gamma)}\int_0^{ s} N(v(u))\de u.
  $$
  The map $\tau$ provides the required reparametrization  by defining the horizontal curve $\widetilde\gamma:[0,T]\to \Hn$
  with $\widetilde\gamma(\sigma)=\gamma( \tau^{-1} (\sigma))$. To see this notice that $\widetilde\gamma(0)=e,\ \widetilde\gamma(T)=g$ and $\gamma( s)=\widetilde\gamma(\tau(s))$.

  Differentiating this relation with respect to $s$ we obtain:

 $$ \frac{d \gamma_I(s)}{d s} = \frac{d\tilde \gamma_I}{d\sigma}(\tau(s))\cdot \tau'(s)=
\frac{d \tilde \gamma_I}{d\sigma}(\tau(s))\frac{T}{\ell(\gamma)} N(v(s),$$
for a.e. $s\in [0,T]$.
Applying the function $N$ to both side of this equality and reordering we obtain:
$$ N(\tilde v(\sigma)) = N\left( \frac{d\tilde \gamma_I}{d\sigma}(\sigma)\right) = \frac{\ell(\gamma)}{T}$$ for a.e. $\sigma \in [0,T]$.

  \QED

Now, let us consider the problem
  \begin{equation}\label{geodesics square}
\left\{
  \begin{array}{ll}
     \displaystyle \inf_v \int_0^T F_N(v(s))\de s\\
\dot x_i= -v_i,\qquad i=1,\ldots, n\\
\dot y_i= -v_{n+i},\qquad i=1,\ldots, n\\
\displaystyle \dot{t}=-2\sum_{i=1}^n\left(y_{i}v_i-x_iv_{n+i}\right)\\
\gamma(0)=e\\
\gamma(T)=g
  \end{array}
\right.
 \end{equation}
where $F_N$ is defined by $N$ via \eqref{FN}.

The following proposition states that the two optimal control problems
  \eqref{control problem} and  \eqref{geodesics square} are in fact equivalent:

 \begin{proposition}\label{equivalence} Let $N:\R^{2n}\to \R_{\ge 0}$ be a norm.
 The problems \eqref{control problem} and  \eqref{geodesics square} are equivalent; more precisely, the control $v^*$ is optimal
 for \eqref{control problem} if and only if $v^*$ is optimal for  \eqref{geodesics square}.
   \end{proposition}

 \noindent \textbf{Proof}: Without loss of generality, let us assume $T=1$.
First, let us suppose that $v^*$ is optimal for \eqref{control problem}, i.e.
$$
\int_0^1 N(v^*)\de s\le   \int_0^1 N(v)\de s,
$$
for every admissible control $v$ . By Proposition \ref{reparametrization}, we are in the position to assume that $N(v^*(s))$ is constant. Hence, by Cauchy--Schwartz inequality, we obtain
$$
\int_0^1 (N(v^*))^2\de s= \left(\int_0^1 N(v^*)\de s\right)^2\le  \left( \int_0^1 N(v)\de s\right)^2\le  \int_0^1 (N(v))^2\de s.
$$
This shows that $v^*$ is optimal for the problem \eqref{geodesics square}.

Conversely, let us suppose that $v^*$ is optimal for \eqref{geodesics square} and by contradiction let us assume that there exists an admissible control $\widetilde v$ such that
$$
\int_0^1 N(\widetilde v)\de s<   \int_0^1 N(v^*)\de s.
$$
Again by Proposition \ref{reparametrization}, we may assume  $N(\widetilde v(s))$ constant. The previous inequality and Cauchy--Schwartz inequality give
$$
\int_0^1 (N(\widetilde v))^2\de s= \left(\int_0^1 N(\widetilde v)\de s\right)^2<  \left( \int_0^1 N(v^*)\de s\right)^2\le
\int_0^1 (N(v^*))^2\de s;
$$
this contradicts the optimality of $v^*$ for the problem \eqref{geodesics square}, proving the claim.  \QED

 Let us remark that, as in the classical case where $N$ is the Euclidean norm, there exists an optimal control for the two problems \eqref{control problem} and  \eqref{geodesics square}. To be precise (see
 for example Theorem 4.1 in \cite{FlRi1975})
  \begin{remark}
  Let $N$ be a norm.
The convexity of $N$ and the superlinearity of $F_N$ guarantee the existence of the optimal control for  \eqref{geodesics square}.
\end{remark}

\subsection{The Pontryagin Maximum Principle with  a general norm $N$.}

In what follows we shall focus our attention to the  study of the problem \eqref{geodesics square}. Our approach is based on Pontryagin's Maximum Principle.
In order to do that, we shall introduce the Hamiltonian:
$$ H: \R^{2n+1}\times\R^{2n+2} \times \R^{2n} \to \R, \quad (x,\lambda,v) \to H(x,\lambda,v)$$
as usual by
\begin{eqnarray*}
H(x, \lambda, v) &=&\lambda_0F_N(v)-\sum_{i=1}^{2n}\lambda_iv_i-2 \lambda_{2n+1} \sum_{i=1}^n\left(y_{i}v_i-x_iv_{n+i}\right)
\end{eqnarray*}
where $x=(x_1,\ldots,x_n,y_1,\ldots,y_n,t)$, $\lambda= (\lambda_0,\lambda_1,\ldots,\lambda_{2n+1})$ and
$v= (v_1, \ldots, v_{2n})$.
Let us reorganize the Hamiltonian function by
\begin{eqnarray*}
H(x, \lambda, v)
&=&\lambda_0F_N(v)-\sum_{i=1}^{n}\left[(\lambda_i+2y_i \lambda_{2n+1}) v_i+ (\lambda_{n+i}-2x_i \lambda_{2n+1}) v_{n+i}\right]\\
&=&\lambda_0F_N(v)-a\cdot v,
\end{eqnarray*}
where the function $a:[0,T]\to \R^{2n}$ is given by
\begin{equation}\label{introduce a}
a=\left[\lambda_1+2y_1 \lambda_{2n+1},\ldots ,\lambda_n+2y_n \lambda_{2n+1},\lambda_{n+1}-2x_1 \lambda_{2n+1},\ldots,\lambda_{2n}-2x_n \lambda_{2n+1}\right].
\end{equation}
The Pontryagin Maximum Principle gives a necessary condition for a control $v$ to be optimal for a control problem. Let us only recall that by
giving a control problem with a cost functional
$$
J(v)=\int_0^T L(s,w(s),v(s))\de s
$$
where $w$ is the unique trajectory associated to the admissible control $v$ via the dynamics and the initial/final points,
the classical assumptions on such cost functional $J$  require that the running cost function $L$ is continuous in $(s,w,v)$, differentiable in $w$ for every fixed $(s,v)$, and the derivatives $\frac{\partial L}{\partial w}$ and $\frac{\partial L}{\partial t}$ are continuous as function of all variables (see for example \cite[Theorem 2.2.1]{LeSc2012} or \cite[Theorem 12.10]{AgSa2004}).
In our case the dynamics is given by the horizontality condition for the curve/trajectory and the running cost function is $L(s,w,v)=F_N(v).$ Since $F_N$ is a convex function, clearly it is continuous in $v$ and we need no other assumptions in order to apply the  Pontryagin Maximum Principle. This is the reason that allows us to study the problem \eqref{geodesics square} without additional regularity assumptions on $N$.

 Pontryagin's  Maximum Principle applied to our situation  (see the mentioned books above) gives the following statement:
\begin{thm}\label{pontryagin}
Let $N:\R^{2n} \to \R_{\ge 0}$ be a norm and let us consider the problem \eqref{geodesics square}.
If $v$ is an optimal control, then
 there exists a multiplier $(\lambda_0^v,\lambda^v)\not=(0,0)$ where
 \begin{itemize}
\item[*]  $\lambda_0^v=\lambda_0$ is a constant in $\{0,1\}$
\item[*]  $\lambda^v=\lambda=(\lambda_1,\ldots,\lambda_{2n+1}):[0,T]\to\R^{2n+1}$ is
  an absolutely continuous function
 \item[*]  $(\lambda_0^v,\lambda^v)\not=(0,0)$
 \end{itemize}
  such that, for $s\in[0,T]$ the following properties hold:
\begin{eqnarray}
&&v(s)\in\arg\min_{u\in\R^{2n}} \left(\lambda_0F_N(u)-a(s)\cdot u\right)\label{pontry1}\\
&&\dot\lambda_{2n+1}(s)=\frac{\partial H}{\partial x_{2n+1}}=0\quad\Rightarrow\lambda_{2n+1}=k\label{pontry5}\\
&&\dot\lambda_i(s)=-\frac{\partial H}{\partial x_i}=-2 k v_{n+i}(s),\qquad i=1,\ldots, n
\label{pontry3}\\
&&\dot\lambda_{n+i}(s)=-\frac{\partial H}{\partial x_{n+i}}=2 k v_i(s),\qquad i=1,\ldots, n
\label{pontry4}\\
&&s\mapsto \lambda_0 F_{N}(v(s))-a(s)\cdot v(s)=c\label{pontry6}
\end{eqnarray}
where $k$ and $c$ are constants, and $a$ is as in \eqref{introduce a}.
\end{thm}

 Relation \eqref{pontry6} holds since the problem is autonomous.
Clearly \eqref{pontry3} and \eqref{pontry4}  imply
\begin{equation}\label{lambdas}
\lambda_i(s)=\lambda_i(0)+2 k y_i(s),\qquad \lambda_{n+i}(s)=\lambda_{n+i}(0)-2kx_i(s),\qquad i=1,\ldots, n
\end{equation}
and hence
\begin{equation}\label{a(s)}
a_i(s)=\lambda_i(0)+4ky_i(s),\qquad a_{n+i}(s)=\lambda_{n+i}(0)-4kx_i(s),\quad  i=1,\ldots n.
\end{equation}
In particular we obtain
\begin{equation}\label{a(0)}
a_i(0)=\lambda_i(0),\qquad a_{n+i}(s)=\lambda_{n+i}(0),\quad  i=1,\ldots n.
\end{equation}

We would like to emphasize the following:
\begin{remark}
For every optimal control $v$, the multiplier $(\lambda_0^v,\lambda^v)$ associated to $v$, in general, it is not unique.
By construction, the function $a$ depends on the choice of the multiplier $(\lambda_0^v,\lambda^v)$. Such function $a$ is absolutely continuous.
\end{remark}
Let us start our investigation. First we study the case of a general norm and later on we add further assumptions of strict convexity and smoothness.

For our optimal control  we have the following normality property:

\begin{proposition}\label{normality}
Assume that $N:\R^{2n}\to \R_{\ge 0}$ is a norm. Let $v$ be an optimal control for the problem \eqref{geodesics square}. Then $v$ is a normal control, i.e.  $\lambda_0^v =\lambda_0=1$.
\end{proposition}
\noindent \textbf{Proof}:
Let us assume by contradiction that $\lambda_0=0$.   Since $v$ is optimal, the Maximum Principle \eqref{pontry1} guarantees that the $\min$ exists, for every $s\in [0, T]$.

On the other hand,  the function $u\mapsto a(s)\cdot u$ is affine. This implies that the above $\min$ exists only if $a=0$ in $[0,T]$. But then \eqref{a(s)} implies that  $v=0$ in $[0,T]$. This gives  $x_i(0)=x_i(s)$ and $y_i(0)=y_i(s)$
for every $s\in[0,T]$. Clearly this implies that $t=0$ as well. This gives a contradiction to $\gamma(0)\not=\gamma(T)$.
Hence, we conclude that  $\lambda_0=1$. \QED

The following result follows essentially from Proposition \ref{Propo first} and will prove to be useful in the sequel.

\begin{proposition}\label{boh}
Assume that $N:\R^{2n} \to \R_{\ge 0}$ is a norm. Let $v$ be an optimal control for the problem \eqref{geodesics square}.
Then
\begin{itemize}
\item[a.] the Maximum Principle \eqref{pontry1} is equivalent to
\begin{equation}\label{inclusion ok}
a(s)\in \partial F_N(v(s))=N(v(s))\partial N(v(s)),\qquad \forall s\in[0,T];
\end{equation}
which is equivalent to
\begin{equation}\label{inclusion 2 ok}
v(s)\in \partial (F_N)^*(a(s))
= N_*(a(s))\partial N_*(a(s))
,\qquad \forall s\in[0,T].
\end{equation}

\item[b.] there exists a unique constant $R(v)=R>0$ such that
\begin{equation}\label{R 3}
N(v(s))=N_*(a(s))=R,\qquad \forall s\in[0,T].
\end{equation}
\end{itemize}
\end{proposition}
\noindent \textbf{Proof}:
Since $v$ is optimal, we have $\lambda_0=1$ and the Maximum Principle \eqref{pontry1} guarantees that the $\min$ exists, for every $s$.
Using the definition of subdifferential, it is easy to see that a point $v(s)$ realizes the $\min$  in \eqref{pontry1} if and only if \eqref{inclusion ok} holds. The equality in \eqref{inclusion ok} follows from \eqref{chain rule}.

On the other hand, by \eqref{Fenchel}, \eqref{inclusion ok} is equivalent to
$$
v(s)\in \partial (F_N)^*(a(s))
,\qquad \forall s\in [0, T].
$$
and again by \eqref{chain rule}  we obtain \eqref{inclusion 2 ok}.
This proves a.

The convexity of $F_N$ implies, by \eqref{Fenchel}, that \eqref{inclusion ok} is equivalent to
\begin{equation}\label{Fenchel R}
 F_N(v(s))+(F_N)^*(a(s))=v(s)\cdot a(s),\qquad \forall s\in [0,T].
\end{equation}
Note that Proposition \ref{Propo first} gives $(F_N)^*(a(s))=\frac{1}{2}(N_*(a(s)))^2$.
Now  \eqref{pontry6} gives that the function $s\mapsto F_{N_*}(a(s))$ is in fact constant, i.e.
\begin{equation}\label{inclusion 3}
N_*(a(s))=R, \qquad \forall s\in [0,T],
\end{equation}
for some fixed $R\ge 0$. If $R=0$, relation \eqref{inclusion 2 ok} implies $v(s)=0$ in $[0,T]$ contradicting the fact that  $\gamma(T)=0\not=g.$ Thus we have that $R>0$.
\medskip

Now  \eqref{inclusion 2 ok} becomes
$
v(s)\in  R\partial N_*(a(s))
,$ for all $s\in [0,T]$.
This implies $\frac{v(s)}{R}\in \partial N_*(a(s))$ and taking into account \eqref{Eulero} we have
$\frac{v(s)}{R}\cdot a(s)=N_*(a(s)),$
i.e.
\begin{equation}\label{R 2}
v(s)\cdot a(s)=R^2.
\end{equation}
Equation \eqref{Fenchel R} becomes
$$\frac{1}{2}(N(v(s))^2+\frac{1}{2}R^2=R^2,\qquad \forall s.$$
Clearly $R$ depends only of $v$. This proves b.\QED

It is clear that \eqref{R 3} is equivalent to
$$
v(s)\in\partial B_N(0,R),\qquad a(s)\in\partial B_{N_*}(0,R),\qquad \forall s\in[0,T].
$$

\medskip
In what follows we would like to clarify the relation between the solutions of  \eqref{geodesics square} and the concept of geodesics as defined in the introduction of this paper.
In our problem \eqref{geodesics square} we fix $T>0$, the initial point $e$ and a final point $g\in\H^n$.
Let us consider an optimal control $v$ and the associated trajectory $\gamma:[0,T]\to\H^n$ such that $v=-\dot{\gamma}_I$,
$\gamma(0)=e$ and $\gamma(T)=g$. The Pontryagin Maximum Principle   gives us the multiplier $(1,\lambda^v)$ and hence a function $a$.  The tern $(v,\lambda^v,a)$ is in general not unique. It  satisfies
relations \eqref{pontry1}--\eqref{pontry6} and the previous Proposition \ref{boh}, for the same $R>0$. This implies that
$$d_{SF}(0,g)=\int_0^T N(v(s))\de s=TR,$$
where $d_{SF}$ is, as in \eqref{sub-finsler},  the sub--Finsler distance from the points $e$ and $g$.

\noindent
Now we can distinguish  two cases:
\begin{enumerate}
\item if $R=1$, then $\gamma$ is a geodesic according to \eqref{def geodesic}.

\item if $R\not=1$, then $\gamma$ is not a geodesic according to \eqref{def geodesic}. In this case we have to change the parametrization of $\gamma$ to obtain a geodesic. In order to do that we define
$\widetilde\gamma:[0,TR]\to\H^n$ by
\begin{equation}\label{reparametrization bis}\widetilde\gamma(s)=\gamma(s/R).\end{equation}
Clearly ${\rm Im}(\widetilde\gamma)={\rm Im}(\gamma)$, $\widetilde\gamma(RT)=\gamma(T)=g$ and $\widetilde v(s)= v(s/R)/R.$
It is easy to see that  such $\widetilde v$ is an optimal control for  the problem in \eqref{geodesics square}, with the same final point $g$ and final time $RT$.
Since $N$ is homogeneous we obtain, for every $s\in[0,T]$
$$R=N(v(s))=N(R\widetilde v(sR))=R N(\widetilde v(sR)).$$
Hence the \lq\lq$R$\rq\rq\ associated to $\widetilde v$  is 1. Now we obtain
$$d_{SF}(0,g)=\int_0^{TR} N(\widetilde v(s))\de s=TR.$$
The new curve $\widetilde\gamma$ is really the geodesic in the sense of \eqref{def geodesic}.
\end{enumerate}
Conversely, if we have a geodesic $\gamma:[0,T]\to\H^n$ from $0$ to $g$, then $v=-\dot\gamma_I$ is an optimal control for \eqref{geodesics square} and the uniqueness of $R$ in Proposition \ref{boh} implies easily that $R=1$.

We would like to formulate the above observation in the following:
\begin{remark}\label{relation optimal control and geodesic}
A horizontal curve $\gamma:[0,T]\to\H^n$, with initial point the origin, is a geodesic if and only if the associated $v=-\dot\gamma_I$ is a optimal control for \eqref{control problem} with $\gamma(T)=g\not= e$ and $R=1$.

 \noindent
Moreover, if for a horizontal curve $\gamma:[0,T]\to\H^n$, with initial point the origin, the associated $v=-\dot\gamma_I$ is an optimal control for \eqref{geodesics square} with $\gamma(T)=g\not=e$ and $R\not=1$, after a reparametrization of $\gamma$ (as in \eqref{reparametrization bis}) we obtain a geodesic.
\end{remark}
In the sequel, when we deal with a finite geodesic $\gamma$, then $N(-\dot{\gamma}_I(s))=R=1$.
Furthermore, in the case when $\gamma:[0,\infty)\to\H^n$ is an infinite geodesic, then for every fixed $T>0$ the curve $\gamma\bigl|_{[0,T]}$ is a finite geodesic, which implies that we have that $N(-\dot{\gamma}_I(s))=R=1$ for all $s\in [0, \infty)$.



In the following we shall add the assumption of strict convexity of the norm $N$.
\subsection{$N$ is a strictly convex norm.}

Let recall  that the strict convexity of $N$ implies that $N_*$ is smooth.

\begin{proposition}\label{boh bis}
Assume that $N:\R^{2n}\to \R_{\ge 0}$ is a strictly convex norm. Let $v$ be an optimal control for the problem \eqref{geodesics square}.
Then
\begin{itemize}

\item[c.] for every $a(s)$ there exists a unique $v(s)$ such that inclusion  \eqref{inclusion 2 ok} holds in the form of an equality; more precisely we obtain
\begin{equation}\label{v via nabla}
v(s)=R\nabla N_*(a(s)),
\end{equation}
where $R>0$ is the constant from \eqref{R 3}.
\item[d.] If $k=0$, then the unique solutions $\gamma=(x_1,\ldots,x_n,y_1,\ldots,y_n,t)$ of  \eqref{pontry1}--\eqref{pontry6} are horizontal segments.

\end{itemize}

\end{proposition}
\noindent \textbf{Proof}:
The strict convexity of $N$ as a norm implies the strict convexity of the function $F_N$ (see Proposition \ref{Propo zero bis}): this gives that for every $a(s)\in \R^{2n}$ the value of
$$\min_{u\in\R^{2n}}\left( F_N(u)-a(s)\cdot u \right)$$
is achieved at a unique point $v(s) \in \R^{2n}$. This proves c.

Let $k=0$.  It is easy to see, by \eqref{a(s)}, that we obtain $a_i(s)=\lambda_i(0)$ for $s\in[0,T],\ i=1,\ldots 2n$.
 In the Maximum Principle \eqref{pontry1} the $\min$ exists, is unique and since $a(s)$ does not depend on $s$,
we obtain for such min that $v=v(s)$ does not depend on $s$. Hence we have that $\gamma$ is a horizontal segment. \QED

This argument does not work without the strict convexity assumption on $N$: more precisely if $k=0$ and if $N$ and thus $F_N$ is only a convex function, it is not always possible to guarantee that for every
fixed $a(s)$, the function
$$u\mapsto F_N(u)-a(s)\cdot u$$
has a unique minimum, for every fixed $s$.

 Let us recall that in general the subdifferential as a set-valued mapping $a\mapsto\partial N_*(a)$ is upper continuous. In our case the subdifferential is single valued $\partial N_*= \nabla N_*$ and consequently we obtain that $a\to \nabla N_*(a)$ is continuous. Furthermore, using the fact that the function $s\to a(s)$ is  absolutely continuous, we obtain by \eqref{v via nabla}, the following:
\begin{remark}\label{v C zero}
If $N:\R^{2n}\to \R_{\ge 0}$ is a strictly convex norm and $v$ is an optimal control for the problem \eqref{geodesics square}, then $v$ is continuous.
\end{remark}

\section{Proof of the main theorem}

In this section we give the proof of Theorem \ref{main}. For the sake of a better understanding we have decided to give the proof in two stages. In the first stage we consider the special case when $N$ is a smooth norm. In this case the proof is easier as we can show uniqueness of the multiplier $\l$ coming from Pontryagin Maximum Principle.

The general case, when $N$ is strictly convex but possibly non-smooth, some technical complications arise as in general the multiplier is not unique. In the second step we show how to overcome this difficulty and prove Theorem \ref{main} also in this case.

\subsection{$N$ is a strictly convex and smooth norm.}

\begin{proposition}\label{prop smooth}
Assume that $N: \R^{2n}\to \R_{\ge 0}$ is a strictly convex and smooth norm. Let $v:[0,T]\to \R^{2n}$ be an optimal control for the problem \eqref{geodesics square}.
Then we have that
\begin{itemize}
\item[e.] for every $v(s)$ there exists a unique $a(s)$ such that inclusion \eqref{inclusion ok} holds;

\item[f.] the associated multiplier $\lambda^v$ (and hence the function $a$ in \eqref{a(s)}) is unique.
\end{itemize}
\end{proposition}
\noindent \textbf{Proof}:
Since $N$ is differentiable except of  the origin, we know that $\partial N(v(s))= \{\nabla N (v(s))\}$  for every $s\in [0,T]$. More precisely \eqref{inclusion ok} becomes the equality
$$a(s)= R\nabla N(v(s)),\qquad \forall \ s\in [0,T].$$ Hence e. is proved.

Let $v:[0,T]\to \R^{2n}$ be the optimal control of the problem  \eqref{geodesics square} and let us assume that there exist two multipliers $\lambda^v,\ \widetilde\lambda^v:[0,T]\to\R^{2n+1}$ associated to $v$.
Let $a$ and $\widetilde a$ defined by \eqref{a(s)} by  $\lambda^v$ and $\widetilde\lambda^v$ respectively. Now e. implies
\begin{equation}\label{equal a}
a(s)=\widetilde a(s),\qquad s\in[0,T]
\end{equation}
and hence, $\dot a(s)=\dot{\widetilde a}(s)$. By  \eqref{a(s)} this implies that  $4\lambda_{2n+1}\dot y_i=4\widetilde{\lambda}_{2n+1}\dot y_i$ and $4\lambda_{2n+1}\dot x_i=4\widetilde{\lambda}_{2n+1}\dot x_i$.
Since the associated trajectory is not a constant curve, there exists a value $s\in [0,T]$ and an index $i$ with the property that $\dot x_i(s)$ or $\dot y_i(s)$ does not vanish. This implies $k=\lambda_{2n+1}=\widetilde{\lambda}_{2n+1}$.
In \eqref{equal a} we obtain
$$
\lambda_i(0)+4ky_i(s)=\widetilde \lambda_i(0)+4ky_i(s),\qquad
\lambda_{n+i}(0)-4kx_i(s)=\widetilde \lambda_{n+i}(0)-4kx_i(s),
\qquad i=1,\ldots,n.
$$
Relation \eqref{lambdas} gives f.\QED

As we will show, if $N$ is only a strictly convex norm, then Proposition \ref{prop smooth} fails (see Example \ref{example 2}).

In what follows we shall consider infinite geodesics $\gamma: \R \to \Hn$. We restrict $\gamma$ to the half-line $[0, \infty)$ and we call this restriction also an infinite geodesic. All our considerations can be repeated to the restriction of $\gamma$ to negative parameter values and so it is enough to consider only the restriction $\gamma: [0,\infty) \to \Hn$.
\medskip

As a next step we shall restrict $\gamma$ to the finite interval $[0,T]$ for some $T>0$ we can conclude our curve is a solution to the optimal control problem \eqref{geodesics square}. The main technical difficulty in our analysis  is the possible dependence on the value of $T$
of the multiplier $\lambda^v$ from Pontryagin's theorem. Indeed, Pontryagin's theorem does not guarantee the uniqueness of the multiplier $\lambda^v$ for a given optimal control $v$ and it can happen that for two different values of say $T'>T$ we obtain two different multipliers $\lambda^{v,T}$ and $\lambda^{v, T'}$ that do not necessarily coincide on the the interval $[0,T]$. The following proposition shows how to deal with this problem if the norm $N$ is smooth.

Let $\gamma:[0,\infty)\to\H^n$ be a infinite geodesic. We say  that $\lambda$ is an  \emph{infinite multiplier}  if for every fixed $T>0$ the function $\lambda\bigl|_{[0,T]}$ is a multiplier associated to the optimal control $-\dot{\gamma}_I\bigl|_{[0,T]}$ for the problem \eqref{geodesics square} with $\gamma(T)=g$.

\begin{proposition}\label{infinite}
Let $N: \R^{2n}\to \R_{\ge 0}$ be a strictly convex and smooth norm. Let $\gamma:[0,\infty)\to\H^n$ be a infinite geodesic sub-Finsler geodesic with respect to the sub-Finsler metric associated to $N$. Then there exists a unique infinite multiplier
$\lambda:[0,\infty)\to\R^{2n+1}$, and hence a unique function $a:[0,\infty)\to\R^{2n}$ via \eqref{a(s)} such that $N_*(a(s))=1$.
\end{proposition}

\noindent \textbf{Proof}: Let $\gamma: [0,\infty)\to\H^n$ be an infinite sub-Finsler geodesic associated to $N$ by \eqref{def geodesic} and \eqref{sub-finsler}.


For some positive $T>0$, the function $\gamma\bigl|_{[0,T]}$ is a geodesic from $\gamma(0)=e$ to $\gamma(T)$, with horizontal velocity $v\bigl|_{[0,T]}$. Since $N$ is smooth (see Proposition \ref{prop smooth}), there exists a unique pair $(\lambda^T,a^T)$ associated to such horizontal velocity such that
$$
a^T(s)=\nabla N(v(s)),\qquad s\in [0,T].
$$

Notice that if we take $T'>T$ then the above relation will also be satisfied by $a^{T'}$ on $[0,T]$ i.e.
$$
a^{T'}(s)=\nabla N(v(s)),\qquad s\in [0,T].
$$

From the above relation we conclude that $a^T= a^{T'}$ on the common interval $[0,T]$ showing that in fact $a^T$ does not really depend on $T$ and we are in the position to define the a $a$ in $[0,\infty)$ such that
$$ a(s) =  \nabla N(v(s)),\qquad s\in [0,\infty).$$
Uniqueness of $a$ can be used to conclude the uniqueness of the multiplier $\l$ as follows. Notice first that relation \eqref{a(s)} gives
$$
a_i(s)=\lambda_i^T(0)+4\lambda^T_{2n+1}y_i(s),\qquad a_{n+i}(s)=\lambda_{n+i}^T(0)-4\lambda^T_{2n+1}x_i(s),\d  i=1,\ldots n,\qquad s\in[0,T].
$$
Taking derivatives w.r.t. $s$  we obtain
$$ \dot a(s) = 4 \lambda^T_{2n+1} J_n v(s), \qquad s\in[0,T],$$
where $J_n$ is the standard symplectic matrix in $\R^{2n}$. Since there exists a value $s_0 \in (0, \infty)$ such that $v(s_0) \not = 0$ we conclude that  $\lambda^T_{2n+1}$ does not depend on $T$ for $T>s_0$. Setting  $\lambda^T_{2n+1}= k$ in the previous relation, it is easy to see the $\lambda^T$ does not depend on $T$ for $T>s_0$ proving the claim.
\QED

The following statement is a special case of Theorem \ref{main} under the additional assumption that $N$ is smooth.

\begin{proposition}\label{final}
Let $N :\R^{2n}\to \R_{\ge 0}$ be a strictly convex and smooth norm. Let
$$\gamma=(x_1,\ldots,x_n,y_1,\ldots,y_n,t):\R\to\H^n$$
 be an  infinite geodesic with respect to the associated sub-Finsler metric. Then $\gamma$ is a horizontal line.
\end{proposition}
\noindent \textbf{Proof}: We consider again the restriction $\gamma:[0, \infty) \to \Hn$.
  Proposition \ref{infinite} implies that we can associate a unique function  $a$ satisfying \eqref{a(s)} with $N_*(a(s))=1$ for $s\ge 0$.

  From Proposition \ref{infinite} we can associate a values $\lambda_i(0)$ for $i=i, \ldots 2n$ and $\lambda_{2n+1}= k$ satisfying
 \begin{equation} \label{a-gamma}
a_i(s)=\lambda_i(0)+4ky_i(s),\quad a_{n+i}(s)=\lambda_{n+i}(0)-4kx_i(s),\qquad  i=1,\ldots n,\quad s\in[0,\infty).
\end{equation}
There are two cases to consider: $k=0$ and $k\not =0$.

 In the first case we obtain that $a$ is
constant vector. On the other hand we have the equation $ a= a(s) =  \nabla N(v(s)), s\in [0,\infty)$. Using this relation we obtain that $v$ must also be a constant. This implies that
$\{ \gamma(s): s\in [0, \infty)\}$ is a horizontal half-line. The same argument can be done for negative values of $s$ as well and we can conclude that $\gamma$ is a horizontal line.

\medskip

In what follows we shall prove that the other possibility: $k \not=0$ will lead to a contradiction.
To see this we use first that $N_*(a(s))=1$ and therefore $a$ takes its values in a compact set.
Then, the relation \eqref{a-gamma} can be used to express $\gamma$ in terms of $a$.
More precisely we obtain:
\begin{equation} \label{a-gamma1}
\frac{a_i(s)-\lambda_i(0)}{4k}=y_i(s),\quad \frac{a_{n+i}(s)-\lambda_{n+i}(0)}{-4k}=x_i(s),\qquad   i=1,\ldots n,\quad s\in[0,\infty).
\end{equation}
 This implies that there exists a constant $C>0$ such that
  \begin{equation}\label{bounded}
  \| (x_1(s), \ldots, x_n(s), y_1(s), \ldots, y_n(s)\| = \|(\gamma_1(s),\ldots,\gamma_{2n}(s))\|\le C,\qquad \forall s\ge 0.
\end{equation}
  This means that the restriction of our geodesic $\gamma$ to $[0, \infty)$ lies in infinite cylinder of a fixed radius $C >0$.
 Let us define $\{\gamma_k\}_{k\ge 1}$, $\gamma_k:[0,\infty ) \to\H^n$ by
$$
\gamma_k(s)=\delta_{\frac{1}{k}}(\gamma(sk)),\qquad s\in [0, \infty).
$$
In what follows we shall check  that $\gamma_k:[0,\infty)\to \Hn$ is a infinite geodesic for every $k$.

In order to do that observe first that since $\gamma:[0,\infty)\to \Hn$ is a geodesic with respect to $d_{SR}$ we have
$$ d_{SF}(\gamma(s), \gamma(s')) = |s-s'| , \qquad \text{for all} \quad s, s' \in [0, \infty).$$
Using the homogeneity of the metric $d_{SF}$ we can infer that
\begin{eqnarray*}
d_{SF}(\gamma_k(s), \gamma_k(s')) = d_{SF}(\delta_{\frac{1}{k}}\gamma (ks), \delta_{\frac{1}{k}}\gamma (ks'))
= \frac{1}{k}d_{SF}(\gamma(ks), \gamma(ks'))=\\ =\frac{1}{k}|ks-ks'| = |s-s'| \quad \text{for all} \quad s, s' \in [0, \infty).
\end{eqnarray*}

This relation shows that $\gamma_k:[0,\infty)\to \Hn$ is an infinite geodesic for every $k$ as claimed. At the same time this also shows that the family of functions $\{\gamma_k\}_{k\ge 1}$, $\gamma_k:[0,1] \to\H^n$ is uniformly bounded and equicontinuous. Uniform boundedness of the family follows from the fact that  $\gamma_k(0)=0$ for all $k$. Applying the theorem of Arzel\'a-Ascoli we obtain a subsequence $\gamma_{n_k}:[0,1] \to \Hn$ converging uniformly to a function $\widehat\gamma:[0,1] \to \Hn$
Using the fact that
$$d_{SF}(\gamma_{n_k} (s),\gamma_{n_k} (s')) = |s-s'| \quad \text{for all} \quad s, s' \in [0, 1],$$
and taking the limit in this equation as $k \to \infty$ we obtain
$$d_{SF}( \widehat\gamma(s),\widehat\gamma(s'))= |s-s'| \quad \text{for all} \quad s, s' \in [0, 1].$$
This shows that $\widehat\gamma:[0,1] \to \Hn$ is a gedoesic.

\medskip
  On the other hand, relation \eqref{bounded} implies that
$\widehat\gamma=(\widehat\gamma_1,\ldots,\widehat\gamma_{2n+1})$ will satisfy $\widehat\gamma_i=0$ for $i=1\ldots 2n$. This means that the curve $\widehat\gamma$ is contained in the vertical axis. It is clear also that the the image $\widehat \gamma$ is not reduced to a single point and so it must contain a non-trivial interval of the vertical axis. However, it is well known that the Hausdorff dimension with respect to $d_{SF}$ of any non-trival interval of the vertical axis is equal to $2$. (For this statement and more general related results we refer to \cite{BalTysWar2005}). This clearly   contradicts the fact that  $\widehat\gamma$ is a geodesic.
\QED

\subsection{$N$ is a strictly convex and non-smooth norm.}

In this subsection we show the modifications necessary to prove Theorem \ref{main}  in the general case when the norm $N$ is strictly convex and not necessarily smooth.



 In fact the proof follows along the same lines as the one of Proposition \ref{final}. The only missing piece from the puzzle is the following boundedness property of infinite geodesics.

\begin{proposition}\label{proj-estimate}
Let $N:\R^{2n} \to \R_{\ge 0}$ be a strictly convex and possibly non-smooth norm. Let
$$\gamma=(x_1,\ldots,x_n,y_1,\ldots,y_n,t):\R\to\H^n$$ be an  infinite geodesic with respect to the sub-Finsler metric associated to $N$. Assume in addition that $\gamma$ is not a
horizontal line.
Then there exists a constant $C>0$ such that
\begin{equation}\label{non-smooth-bounded}
  \| (x_1(s), \ldots, x_n(s), y_1(s), \ldots, y_n(s))\| = \|(\gamma_1(s),\ldots,\gamma_{2n}(s))\|\le C,\qquad \forall s\in \R.
\end{equation}

\end{proposition}
\noindent \textbf{Proof}:
We shall prove the estimate of \eqref{proj-estimate} only of values $s\geq 0$ For negative values of $s$ the argument is similar. Let us restrict $\gamma$ to the finite interval $[0,T]$
for some $T>0$.

Let us recall  that  relation \eqref{a(s)} gives
$$
a_i^T(s)=\lambda_i^T(0)+4\lambda^T_{2n+1}y_i(s),\quad a_{n+i}^T(s)=\lambda_{n+i}^T(0)-4\lambda^T_{2n+1}x_i(s),\quad  i=1,\ldots n,\quad s\in[0,T].
$$
Setting $s=0 $ in the above relations and using that $\gamma(0)=e$ we conclude that
$a^T_i(0) = \lambda^T_i(0)$, for every $i$, and so we can write these relations
 in shorthand as
\begin{equation} \label{shorthand-gamma}
a^T(s) = a^{T}(0) + 4 k^T J_n \gamma_I(s), \qquad s\in [0,T].
\end{equation}
where $J_n$ is the standard symplectic matrix in $\R^{2n}$.
Let us choose a parameter value $s_0$ for which $v(s_0) \not = v(0)$ and $\gamma_I(s_0) \not = 0$.  Since $\gamma$ is not a line, such value $s_0$ can be found for a sufficient large $T$. Let us write the above relation \eqref{shorthand-gamma} at $s_0$ as
\begin{equation} \label{shorthand-gamma-0}
a^T(s_0) - a^{T}(0) = 4 k^T J_n \gamma_I(s_0), \qquad s\in [0,T].
\end{equation}
Notice that $a^T(s_0) \in \partial F_N(v(s_0))$, $a^T(0)\in \partial F_N(v(0))$. Furthermore, observe that $\partial F_N(v(s_0))$ and $\partial F_N(v(0))$ are two compact and convex sets. By the strict convexity of $N$ the intersection of $\partial F_N(v(s_0))$ and $\partial F_N(v(0))$ is empty according to Proposition \ref{Propo zero bis}.
This implies that
\begin{equation} \label{min-estim}
c:=\min \{\|u-u'\|: u \in \partial F_N(v(s_0)), u'\in \partial F_N(v(0))\} >0,
\end{equation}
which yields the estimate
\begin{equation} \label{norm-estim}
0<c \leq \|a^T(s_0) - a^{T}(0)\|
\end{equation}
for all values of $T>s_0$.
Taking the norm in \eqref{shorthand-gamma-0} we obtain
$$c \leq 4 |k^T| \| \gamma_I(s_0) \|,$$
 which implies that
 \begin{equation} \label{k-estim}
 |k^T| \geq \frac{c}{4\| \gamma_I(s_0)\|}, \qquad T >s_0.
 \end{equation}

Now let us turn back to the relation \eqref{shorthand-gamma}: using that $N(v(s))  = 1$, $a^T(s)\in \partial F_N(v(s))$ and the fact that the operator $v\to \partial F_N(v)$ is upper semicontinuous (see \cite[Proposition 2.1]{AlAmCa1992}) we obtain that $a^T(s)$ takes values in a fixed compact set that is independent on $T$ and $s$.
We can now write  \eqref{shorthand-gamma} in the form
$$\frac{a^T(s) - a^{T}(0)}{4 k^T}= J_n \gamma_I(s), \qquad s\in [0,T].$$
Taking into consideration the above discussion and the estimate \eqref{k-estim} we can conclude that the quantity on the left hand side of the above relation takes values in a fixed compact set independently on $s$ and $T$ proving our claim.
\QED

\section{Examples and final remarks}
Let $p\in\left[1,\infty\right]$, and let $\|\cdot \|_{p}$ be the $p$-norm on $\mathbb{R}^{2n}$  and $a\in\left(0,\infty\right)$. Then one can check (see \cite{BalFasSob2017}) that the function
\begin{align}
\hat N= N_{p,a}:\mathbb{H}^n\to\mathbb{R},\ (z,t)\mapsto \max\left\{ ||z||_{p},a\sqrt{|t|}\right\},\nonumber
\end{align}
defines a left-invariant  norm on $\mathbb{H}^n$, if either $ 1\le p\le 2\ \text{and}\ 0<a\le1$ or if $ 2<p\le\infty\ \text{and}\ 0<a\le n^{1/p-1/2}$.

In \cite{BalFasSob2017} it was proved that that $\hat N= N_{p,a}$ has the $GLP$ iff $p\in (1,+\infty)$.
This can be obtained also as an immediate corollary of Theorem \ref{main} observing that the norm $N$ associated to $\hat N= N_{p,a}$ is the usual $p$-norm: $N(z)= ||z||_p$. Since $||\cdot||_p$ is strictly convex if and only if $p\in (1,+\infty)$ we obtain this result as an immediate corollary of Theorem \ref{main}.  Notice, that for $p=\infty$ or $p=1$ admits infinite geodesics which are not lines, $\hat N= N_{p,a}$ and in fact there exist in this setting isometric embeddings which are non-linear.

\medskip

The above example is a particular case of a more general phenomena. Observe first that  non-strict convexity of a given norm $N$ in $\R^{2n}$ is equivalent with the existence non-linear of infinite geodesics in the normed space $(\R^{2n}, d_N)$. Next, we notice that  horizontal lifts of  infinite geodesics in $(\R^{2n}, d_N)$ will be infinite geodesics in associated sub-Finsler space $(\H^n, d_{SR})$. This shows that there exists an abundance of  non-linear infinite geodesics that arise naturally in $(\H^n, d_{SR})$ whenever the original norm $N$ is not strict convex. For more details about the above argument we refer to the proof of Proposition 3.14 in \cite{BalFasSob2017}.

\medskip
In the first Heisenberg group $\mathbb{H}^1$, the classification of the geodesics with respect to a sub-Finsler distance associated to a norm $N$ is related to the following \emph{isoperimetric problem} on the Minkowski plane $(\mathbb{R}^2,N)$: given a number $A$ find a closed path through $0$ of minimal $N$-length which encloses (Euclidean) area $A$. To describe the solution, we recall the following notation for the closed unit ball and dual ball in $(\mathbb{R}^2,N)$:
\begin{displaymath}
B:= \{z\in\mathbb{R}^2:\; N(z) \leq 1\}\quad\text{and} \quad B^{\circ}:=\{w:\; \langle w,z\rangle\leq 1:\; z\in B\}.
\end{displaymath}
The \emph{isoperimetrix} $I$ is the boundary of $B^{\circ}$ rotated by $\pi/2$, and it can be parameterized as a closed curve.
Buseman \cite{Bus1947}  proved that the solution to the above stated isoperimetric problem is given by (appropriate dilation and translation) of the isoperimetrix.  Note that if $N$ is strictly convex, then $I$ is of class $\mathcal{C}^1$. Considering the associated sub-Finsler space
$(\mathbb{H}^1, d_{SF})$ we can conclude that horizontal lifts of the isoperimetrix a will be in fact geodesics in $(\mathbb{H}^1, d_{SF})$. Let us note that this classical result of Busemann follows rather easily from the control-theoretical approach, in particular Proposition \ref{boh}. In fact we can observe that according to Proposition \ref{boh} the geodesics in $(\Hn, d_{SR})$
are also curves on the boundary of the dual ball, $B^{\circ}$ where
\begin{displaymath}
B:= \{z\in\mathbb{R}^{2n}:\; N(z) \leq 1\}\quad\text{and} \quad B^{\circ}:=\{w:\; \langle w,z\rangle\leq 1:\; z\in B\}.
\end{displaymath}
To be clear, in the case $n=1$, if we have a geodesic then in Proposition \ref{boh} we obtain $N_*(a(s))=N(v(s))=R=1$; this implies,
 using \eqref{a(s)},
 $$
(\lambda_1(0)+4ky(s),\lambda_{2}(0)-4kx(s)) \in \partial B_{N_*},\qquad\forall  s.$$
If $k\not =0,$ this is equivalent to
$$(y(s), -x(s)) \in \frac{1}{4k}\Bigr(\partial B_{N_*}-(\lambda_1(0), \lambda_2(0))\Bigl) ,\qquad\forall  s.$$
And this is equivalent to the fact that
$(x(s),y(s))$ lies in the set described by the Busseman result,
i.e. considering  the set $\partial B_{N_*}$, a translation, a dilation and a rotation of $\pi/2$.
For $n=1$ the whole boundary is curve and so the geodesics in $(\H^1, d_{SF})$ can be explicitly characterised.

\medskip

In higher dimensions such explicit characterisation cannot be expected for a general norm $N$. In certain particular cases this is still possible. An example of this explicit characterisation is the well-known case of the sub-Riemannian geodesics when $N$ is the usual Euclidean norm $N(z) = ||z||_2$ given by Monti \cite{Mon2000} who also used control theoretical approach. In this case the optimal control problem \eqref{geodesics square} can be explicitly solved and the expressions of the sub-Riemannian geodesics are precisely computed \cite{Mon2000}.

\medskip

 As we have showed, the case of $N$ a strictly convex and non-smooth norm may lead to the non-uniqueness and ambiguity of the multiplier coming from Pontryagin's theorem.
 In the next example we show that this situation actually happens.  We construct a solution of the optimal control problem \eqref{geodesics square} such that there exist multiple choices for $a(s)$ such that inclusion \eqref{inclusion ok} holds, since the multiplier $\lambda^v$ is not unique. From this we can see that
 without the assumption on the smoothness of the norm, Proposition \ref{prop smooth} does not hold.

\begin{example}[strictly convex and non-smooth norm]\label{example 2}
Let $N:\R^2\to\R_{\ge 0}$ defined by ,
$$N(z)=|x|+\sqrt{2x^2+y^2},$$
for $z=(x,y)\in\R^2$.
\end{example}
It is an exercise to show that the shape of $\partial B_N$ is a sort of an American football ball, non smooth in the points $(0,\pm 1)$: more precisely we have
$$\partial B_N=\{(x,y)\in\R^2:\ x^2+y^2+2|x|-1=0\}.$$
 We can calculate the dual norm using its definition via some rather tedious calculation to obtain
\begin{equation}\label{example 1 00}
 N_*(x,y)=\left\{
\begin{array}{ll}
\displaystyle -|x|+\sqrt{2}\sqrt{x^2+y^2}& {\rm if}\ |x|\ge |y| \\
\displaystyle |y|& {\rm if}\ |x|<|y|\\
\end{array}
\right.
\end{equation}
with
$$\partial B_{N_*}=
\Bigl\{(x,y)\in\R^2:\ |x|\le 1,\ |y|=1\Bigr\}\cup \Bigl\{(x,y)\in\R^2:\ |x|\ge 1,\ x^2+2y^2-2|x|=1\Bigr\}$$
and
\begin{equation}\label{example 1 01}
\nabla N_*(x,y)=\left\{
\begin{array}{ll}
\displaystyle \left(-\texttt{\rm sgn}(x)+\frac{\sqrt{2}x}{\sqrt{x^2+y^2}},\frac{\sqrt{2}y}{\sqrt{x^2+y^2}}\right)& {\rm if}\ |x|\ge |y| ,\\
\displaystyle \left(0,\texttt{\rm sgn}(y)\right)& {\rm if}\ |x|<|y|.\\
\end{array}
\right.
\end{equation}

 We claim that the curve $\gamma:[0,\tau]\to \H^1$  defined by
\begin{equation}\label{example 1 0}
\gamma(s)=\left\{
\begin{array}{ll}
\displaystyle (0,-s,0)& {\rm if}\ s\in[0,1] \\
\displaystyle \left(-1+\sqrt{(1+2\theta(s)-\theta^2(s))/2},-\theta(s),\gamma_3(s)\right)\qquad& {\rm if}\  s\in(1,\tau]\\
\end{array}
\right.
\end{equation}
is a geodesic. The function $\gamma_3$ is defined via the horizontality condition on $\gamma,$ while the point $\tau$ and the function $\theta$ are defined as follows. Consider the Cauchy problem:
\begin{equation}\label{example 1 1}
\left\{
\begin{array}{l}
\displaystyle \dot\theta=
\frac{\sqrt{2+4\theta-2\theta^2}}{1+\theta}  \\
\theta(1)=1.\\
\end{array}
\right.
\end{equation}
The solution is given by the relation
\begin{equation}\label{theta}
s-\sqrt{2}\arcsin\left(\frac{\theta(s)-1}{\sqrt{2}}\right)+\frac{1}{2}\sqrt{2+4\theta(s)-2\theta(s)^2}=2.
\end{equation}
This relation defines a differentiable, increasing function $\theta:[1,\tau]\to \R$, with $\tau=2+\pi/\sqrt{2}$, such that $\theta(\tau)=1+\sqrt{2}.$

The first step in order to show that $\gamma$ is a geodesic is to prove that it is a solution of the Pontryagin system \eqref{pontry1}--\eqref{pontry6}. Equivalently $\gamma$ is a curve such that the triple $(v=-\dot\gamma_I,\lambda,a)$ satisfies \eqref{pontry1}--\eqref{pontry6} for some $a$ and $\lambda$.

Clearly, by \eqref{example 1 0}, we have
\begin{equation}\label{example 1 1 0}
-\dot\gamma_I(s)=v(s)=\left\{
\begin{array}{ll}
\displaystyle (0,1)& {\rm if}\ s\in[0,1] \\
\displaystyle \left(\frac{\dot\theta(\theta-1)}{\sqrt{2+4\theta-2\theta^2}},\dot\theta(s)\right)\qquad& {\rm if}\  s\in(1,\tau] .\\
\end{array}
\right.
\end{equation}
Let us consider the function $\lambda:[0,\tau]\to \R^3$ defined as in \eqref{lambdas} choosing $\lambda_1(0)=0,\ \lambda_2(0)=1$ and $\lambda_3=k=-1/4$, i.e.
$$\lambda(s)=\left(-\frac{1}{2}y(s),1+\frac{1}{2}x(s),-\frac{1}{4}\right).$$
It is immediate to see that for such $\lambda$ conditions \eqref{pontry5}--\eqref{pontry4} are satisfied.

Now, we define the function $a:[0,\tau]\to \R^2$  via \eqref{a(s)}:
\begin{equation}\label{example 1 1 1}
a(s)=\left\{
\begin{array}{ll}
\displaystyle (s,1)& {\rm if}\ s\in[0,1] \\
\displaystyle \left(\theta(s),\sqrt{(1+2\theta(s)-\theta^2(s))/2}\right)& {\rm if}\  s\in(1,\tau].\\
\end{array}
\right.
\end{equation}
We have to prove that relations \eqref{pontry1} and \eqref{pontry6} hold. First, it is easy to see, using the first line in \eqref{example 1 00} and  $\theta(s)>0$, that
\begin{equation}\label{example 1 2}
N_*(a(s))=1,\qquad \forall s\in[0,\tau].
\end{equation}
Now the Maximum Principle  \eqref{pontry1} is equivalent to \eqref{inclusion 2 ok}, i.e. we have to prove that
 $$v(s)\in N_*(a(s))\partial N_*(a(s))=\{\nabla N_*(a(s))\},\qquad s\in[0,\tau],$$
taking into account \eqref{example 1 1 0}, \eqref{example 1 1 1} and \eqref{example 1 01}.
The previous equality for $s\in[0,1]$ is obvious; for $s\in [1,\tau]$ we have to obtain
$$
\left(\frac{\dot\theta(\theta-1)}{\sqrt{2+4\theta-2\theta^2}},\dot\theta(s)\right)
=\left(-1+\frac{2\theta}{1+\theta},\frac{\sqrt{2+4\theta-2\theta^2}}{1+\theta}\right).$$
The two components of this equation are equivalent; together with the condition $\theta(1)=1$, we obtain exactly the Cauchy problem
\eqref{example 1 1}.

Hence  $(v,\lambda,a)$  solves the Pontryagin system \eqref{pontry1}--\eqref{pontry6} and $\gamma$ is candidate to be a geodesic.
Now, from the mentioned Buseman’s theorem \cite{Bus1947}, $\gamma$ is indeed a finite geodesic.

In the spirit of Proposition \ref{proj-estimate}, let us consider $T=1$ and the restriction of $\gamma$ to $[0,T]:$ we have the geodesic
$\gamma\bigl|_{[0,1]}=
\widetilde \gamma$, i.e. the horizontal segment from the origin to the point $(0,-1,0)$. Hence $\widetilde v=-\dot{\widetilde \gamma}=(0,1)$ is a optimal solution for \eqref{geodesics square} with $g=(0,1,0)$. Let us prove that
for every $(\ell,k)$ such that
\begin{equation}\label{example 1 3}
\max\left(|\ell|,\ |\ell-4k|\right)\le 1
\end{equation}
the function $\lambda^{\ell,k}:[0,1]\to \R^3$ defined by
$$\lambda^{\ell,k}(s)=(\ell-2ks,1,k)$$
 is a multiplier associated to the optimal control $\widetilde v$.

 First the, for every $(\ell,k)$, via \eqref{a(s)}, we have the function  $a^{\ell,k}:[0,1]\to \R^2$ defined as
$$a^{\ell,k}(s)=(\ell-4ks,1).$$
 Conditions  \eqref{pontry5}--\eqref{pontry4} are easily satisfied. Condition \eqref{example 1 3} gives that  $N_*(a^{\ell,k}(s))=1$. Hence
 \eqref{inclusion 2 ok} becomes $$\widetilde v(s)=(0,1)=\nabla N_*(a^{\ell,k}(s)),\qquad \forall s\in{0,1}.$$
Using \eqref{example 1 01}, condition \eqref{example 1 3} guarantees the previous equality.
 Hence $\widetilde v$ admits infinitely many multipliers.


\begin{thebibliography}{10}

\bibitem{AgSa2004}
A.~Agrachev and Y.~Sachkov.
\newblock {\em {C}ontrol {T}heory from the {G}eometric {V}iewpoint}.
\newblock Springer, 2004.

\bibitem{ABB}
A. Agrachev, D. Barilari, and U. Boscain,
   \newblock Introduction to geodesics in sub-Riemannian geometry.
   \newblock {\em Geometry, analysis and dynamics on sub-Riemannian manifolds.
      Vol. II}, EMS Ser. Lect. Math., Eur. Math. Soc., Z\"urich, 1--83, 2016.


\bibitem{AlAmCa1992}
G.~Alberti, L.~Ambrosio, and P.~Cannarsa.
\newblock On the singularities of convex functions.
\newblock {\em Manuscripta Math.}, 76:421–435, 1992.

\bibitem{BBDS2017}
D.~Barilari, U.~Boscain, E.~Le Donne, and M.~Sigalotti.
\newblock Sub-Finsler structures from the time-optimal control viewpoint for some nilpotent distributions
\newblock {\em Journal of Dynamical and Control Systems }, 23:3, 547-575, 2017.

\bibitem{BalFasSob2017}
Z.M. Balogh, K.~F$\ddot{\rm a}$ssler, and H.~Sobrino.
\newblock Isometric embeddings into {H}eisenberg groups.
\newblock {\em Geom. Dedicata}, 185, 163--192, 2018.

\bibitem{BalTysWar2009}
Z.M. Balogh, J.T..~Tyson, and B.~Warhurst.
\newblock Sub-{R}iemannian vs. {E}uclidean dimension comparison and fractal geometry on {C}arnot groups.
\newblock {\em Advances in Mathematics}, 220, (2), 560--619, 2009.


\bibitem{BelRis1996}
A.~Bella\"{\i}che and J.-J.~Risler (eds.).
\newblock {\em Sub--{R}iemannian {G}eometry}, volume 144.
\newblock Progress in {M}athematics, Birkh\"auser, 1996.

\bibitem{BeNeOz2003}
D.~P. Bertsekas, A.~Nedic, and E.~Ozdaglar.
\newblock {\em {C}onvex analysis and optimization}.
\newblock Athena Scientific, Belmont, MA, 2003.

\bibitem{Bus1947}
H.~Busemann.
\newblock {T}he {I}soperimetric problem in the {M}inkowski plane.
\newblock {\em American Journal of Mathematics}, 69/4, 1947.

\bibitem{Cl1983}
F.H. Clarke.
\newblock {\em Optimization and Nonsmooth Analysis}.
\newblock John Wiley \& Sons, New York, 1983.

\bibitem{DLR2017}
 E. Le Donne, S. Li, and T. Rajala
\newblock
   {A}hlfors-regular distances on the {H}eisenberg group without
   biLipschitz pieces.
  \newblock {\em Proc. Lond. Math. Soc. (3)},
   115/2, 348--380, 2017,

\bibitem{FlRi1975}
W.H. Fleming and R.W. Rishel.
\newblock {\em Deterministic and {S}tochastic {O}ptimal {C}ontrol}.
\newblock Springer--Verlag, 1975.
	
\bibitem{LeRi2017}
E.~Le Donne and S.~Rigot.
\newblock Besicovitch covering property for homogeneous distances on the
  {H}eisenberg groups.
\newblock {\em J. Eur. Math. Soc.}, 19:1589--1617, 2017.

\bibitem{HaZi2015}
P.~Haj{\l}asz and S.~Zimmerman.
\newblock Geodesics in the {H}eisenberg group.
\newblock {\em Anal. Geom. Metr. Spaces}, 3:325–337, 2015.

\bibitem{LeSc2012}
U.~Ledzewicz and H.~Sch\"attler.
\newblock {\em Geometric Optimal Control}.
\newblock Springer, 2012.

\bibitem{Me1998}
R.E. Megginson.
\newblock {\em An Introduction to Banach Space Theory}.
\newblock Springer, 1988.

\bibitem{Mon2000}
R.~Monti.
\newblock Some properties of {C}arnot--{C}arath\'eodory balls in the
  {H}eisenberg group.
\newblock {\em Rend. Mat. Acc. Lincei}, 11:155--167, 2000.

\bibitem{Pa1989}
P.~Pansu.
\newblock M\'etriques de {C}arnot--{C}arath\'eodory et quasi--isom\'etries des
  espaces sym\'etriques de rang un.
\newblock {\em Ann. of Math.}, 129(2):1--60, 1989.

\bibitem{Pap2005}
A.~Papadopoulos.
\newblock Metric spaces,convexity and non positive curvature
\newblock {\em Lectures in mathematics and theoretical physics, European Mathematical Society}, 2005.

\bibitem{RoWe2004}
R.T. Rockafellar and R.J-B. Wets.
\newblock {\em Variational {A}nalysis}.
\newblock Springer, 2004.

\end{thebibliography}

\end{document}